\documentclass[12pt,letterpaper,reqno]{amsart}

\usepackage[T1]{fontenc}
\usepackage[utf8]{inputenc}
\usepackage{lmodern}
\usepackage{amsmath,amssymb,amsthm,mathtools,mathrsfs}
\usepackage{aliascnt}
\usepackage{microtype}
\usepackage{enumitem}
\usepackage{booktabs}
\usepackage{xcolor}
\usepackage{doi}
\usepackage{hyperref}
\usepackage{tikz}
\usetikzlibrary{positioning,arrows.meta}
\usepackage{bookmark}
\usepackage[capitalize,noabbrev]{cleveref}

\newtheorem{theorem}{Theorem}[section]

\newaliascnt{lemma}{theorem}
\newtheorem{lemma}[lemma]{Lemma}
\aliascntresetthe{lemma}

\newaliascnt{proposition}{theorem}
\newtheorem{proposition}[proposition]{Proposition}
\aliascntresetthe{proposition}

\newaliascnt{corollary}{theorem}
\newtheorem{corollary}[corollary]{Corollary}
\aliascntresetthe{corollary}

\newaliascnt{conjecture}{theorem}
\newtheorem{conjecture}[conjecture]{Conjecture}
\aliascntresetthe{conjecture}

\newaliascnt{problem}{theorem}

\aliascntresetthe{problem}

\theoremstyle{definition}
\newaliascnt{definition}{theorem}

\aliascntresetthe{definition}

\newaliascnt{example}{theorem}

\aliascntresetthe{example}

\newaliascnt{remark}{theorem}

\aliascntresetthe{remark}

\crefname{theorem}{Theorem}{Theorems}
\Crefname{theorem}{Theorem}{Theorems}
\crefname{lemma}{Lemma}{Lemmas}
\Crefname{lemma}{Lemma}{Lemmas}
\crefname{proposition}{Proposition}{Propositions}
\Crefname{proposition}{Proposition}{Propositions}
\crefname{corollary}{Corollary}{Corollaries}
\Crefname{corollary}{Corollary}{Corollaries}
\crefname{conjecture}{Conjecture}{Conjectures}
\Crefname{conjecture}{Conjecture}{Conjectures}

\DeclareMathOperator{\curl}{curl}
\DeclareMathOperator{\dom}{dom}
\DeclareMathOperator{\Tr}{Tr}
\newcommand{\dd}{\,\mathrm{d}}
\newcommand{\C}{\mathbb C}
\newcommand{\R}{\mathbb R}

\begin{document}

\title[Neumann--Dirichlet comparison at the first threshold]
{Neumann--Dirichlet eigenvalue comparison at the first Dirichlet threshold in the plane}

\author[Q.~Tang]{Quanyu Tang}
\address{School of Mathematical Sciences, University of Science and Technology of China, Hefei 230026, P.~R.~China}
\email{tangquanyu827@gmail.com}

\author[H.~Zhang]{Haiqi Zhang}
\address{School of Mathematics, Shandong University, Jinan 250100, P.~R.~China}
\email{ZHQAQ2024@outlook.com}

\begin{abstract}
Let $\Omega\subset\mathbb R^2$ be a bounded connected Lipschitz domain, and let $\{\mu_j(\Omega)\}_{j\geq1}$ and $\{\lambda_j(\Omega)\}_{j\geq1}$ denote the Neumann and Dirichlet Laplacian eigenvalues, respectively, counted with multiplicity. We prove that
$$
\mu_3(\Omega)<\lambda_1(\Omega),
$$
thereby removing the simple-connectivity assumption from the previously known planar result at the first Dirichlet threshold. We also establish a three-spectrum counting inequality relating the Dirichlet, Neumann, and conductivity spectra.
\end{abstract}

% 35P15  	Estimates of eigenvalues in context of PDEs
% 58J50  	Spectral problems; spectral geometry; scattering theory on manifolds [See also 35Pxx]
% 35J05  	Laplace operator, Helmholtz equation (reduced wave equation), Poisson equation [See also 31Axx, 31Bxx]


\maketitle

\section{Introduction}

Let $\Omega\subset\R^d$ be a bounded connected Lipschitz domain. We enumerate
its Neumann and Dirichlet Laplacian eigenvalues, including multiplicities, as
\[
0=\mu_1(\Omega)<\mu_2(\Omega)\leq\mu_3(\Omega)\leq\cdots
\]
and
\[
0<\lambda_1(\Omega)\leq\lambda_2(\Omega)\leq\cdots,
\]
respectively. Comparisons between these two spectra form a classical theme in
spectral geometry. P\'olya proved $\mu_2<\lambda_1$ for planar domains
\cite{Polya1952}. Payne subsequently established $\mu_{k+2}<\lambda_k$ for
sufficiently smooth convex planar domains \cite{Payne1955}. In higher
dimensions, Aviles obtained a curvature-dependent one-step comparison
\cite{Aviles1986}, while Levine and Weinberger proved
$\mu_{k+d}<\lambda_k$ for sufficiently smooth convex domains in $\R^d$ and
observed that approximation yields the corresponding non-strict inequality
for arbitrary bounded convex domains \cite{LevineWeinberger1986}.

A second line of work does not rely on convexity. Friedlander proved
$\mu_{k+1}\leq\lambda_k$ on $C^1$ domains \cite{Friedlander1991}, and Filonov
obtained the strict inequality $\mu_{k+1}<\lambda_k$ for a broad class that
includes bounded Lipschitz domains \cite{Filonov2005}. The
Dirichlet-to-Neumann and Robin-operator approaches were developed further by
Gesztesy and Mitrea and by Arendt and Mazzeo
\cite{GesztesyMitrea2009,ArendtMazzeo2012}, while Safarov developed an abstract formulation in terms of spectral
counting functions, and Safarov--Filonov obtained quantitative asymptotic
estimates for their difference
\cite{Safarov2008,SafarovFilonov2010}. The possibility of retaining
the full dimensional shift without convexity was explicitly discussed by
Benguria, Levitin, and Parnovski \cite{BenguriaLevitinParnovski2009} and is
recorded in \cite[Conjecture~3.2.42]{LevitinMangoubiPolterovich2023}.

The resulting open problem may be stated as follows.

\begin{conjecture}\label{conj:dimensional-shift}
Every bounded connected Lipschitz domain $\Omega\subset\R^d$, $d\geq2$,
satisfies
\begin{equation}\label{eq:dimensional-shift}
\mu_{k+d}(\Omega)\leq\lambda_k(\Omega)
\qquad\text{for every }k\geq1.
\end{equation}
\end{conjecture}

In the plane, \cref{conj:dimensional-shift} is known under a topological
restriction. Rohleder proved \eqref{eq:dimensional-shift} for simply connected planar
Lipschitz domains and obtained strictness whenever the relevant Dirichlet
eigenvalue is simple or the boundary contains a line segment
\cite{Rohleder2025}. In particular, since $\lambda_1$ is simple, his theorem
gives $\mu_3<\lambda_1$ for simply connected planar domains. Fries, Goffeng,
and Miranda subsequently placed the planar argument in the de Rham framework;
here and below we refer to the corrected fourth arXiv version, which removes
an erroneous claim contained in an earlier version \cite{FriesGoffengMiranda2026}. Their
counting formula makes the topological obstruction transparent: the Euler
characteristic contributes positively in the simply connected case but loses
one unit for each hole. Thus the known Hodge-theoretic estimates do not prove
\cref{conj:dimensional-shift} on planar domains of arbitrary connectivity.

Two further Hodge-theoretic developments should be distinguished from the
result below. Muravyev formulated Rohleder's argument through the absolute
Hodge Laplacian and obtained a broader spectral decomposition
\cite{Muravyev2026}. Hua, M\"unch, and Zhang extended the vector-field/Hodge
method to domains in Riemannian surfaces satisfying a curvature condition;
at the first threshold their topology-dependent estimate has the form
$\mu_{3-\beta_1}\leq\lambda_1$, where $\beta_1$ is the first Betti number
\cite{HuaMunchZhang2026}. In the Euclidean plane this recovers the
simply connected estimate when $\beta_1=0$, but it does not give
$\mu_3\leq\lambda_1$ once holes are present. Hence neither result removes the
topological restriction addressed here.

The first Dirichlet threshold nevertheless admits such a comparison, with
strict inequality and with no topological assumption. This is the main result
of the paper.

\begin{theorem}\label{thm:planar-first}
Every bounded connected Lipschitz domain $\Omega\subset\R^2$ satisfies
\begin{equation}\label{eq:planar-first}
\mu_3(\Omega)<\lambda_1(\Omega).
\end{equation}
\end{theorem}

The index in \cref{thm:planar-first} is optimal: already for a disk one has
$\mu_4>\lambda_1$; see \cite{Rohleder2025}. The theorem therefore removes the simple-connectivity
assumption at the first Dirichlet threshold without weakening the conclusion.
The number of Neumann eigenvalues not exceeding $\lambda_1$ has also been studied from
an isoperimetric viewpoint by Cox, MacLachlan, and Steeves and, for convex
domains, by Hatcher \cite{CoxMacLachlanSteeves2019,Hatcher2025}. More broadly,
recent work shows that the natural index gap may grow with $k$ on large classes
of domains \cite{FreitasGama}. In a different, high-dimensional regime,
Filonov obtained index shifts bounded below exponentially in the dimension at
the first Dirichlet threshold on arbitrary domains and at every threshold on
convex domains; in particular, his estimate gives $\mu_3<\lambda_1$ for
$d\geq7$, but it does not cover the planar case \cite{Filonov2026}. Schmatzler has
extended the classical Levine--Weinberger and Friedlander--Filonov mechanisms
to several elliptic operators \cite{Schmatzler2026}.

The proof of \cref{thm:planar-first} is based on a structural comparison with a
third spectrum. Write
$\partial\Omega=\Gamma_0\sqcup\cdots\sqcup\Gamma_g$ for some $g\geq0$, and let
$0<\gamma_1\leq\gamma_2\leq\cdots$ be the eigenvalues associated with the
Dirichlet energy form on zero-mean functions whose traces are constant on
each boundary component. These are
the eigenvalues associated with the conductivity boundary conditions studied
by Bucur, Durus, and Oudet \cite{BucurDurusOudet2008}. The related
zero-mean problem with one globally constant boundary value was studied by
Greco and Lucia \cite{GrecoLucia2008}; the two form domains agree when the
complement of the domain is connected. For $s>0$, define
\[
N_D^{<}(s):=\#\{j:\lambda_j<s\},\qquad
N_{\mathrm c}^{<}(s):=\#\{j:\gamma_j<s\},
\]
and
\[
N_{N,+}^{<}(s):=\#\{j\geq2:\mu_j<s\},\qquad N_N^{<}(s):=\#\{j:\mu_j<s\}=1+N_{N,+}^{<}(s).
\]
Our second main theorem is the following three-spectrum counting inequality.

\begin{theorem}\label{thm:counting-intro}
For every bounded connected Lipschitz domain $\Omega\subset\R^2$ and every
$s>0$,
\begin{equation}\label{eq:counting-intro}
N_{N,+}^{<}(s)
\geq
\max\bigl\{N_{\mathrm c}^{<}(s),\,
2N_D^{<}(s)-N_{\mathrm c}^{<}(s)\bigr\}.
\end{equation}
\end{theorem}

The two alternatives in \eqref{eq:counting-intro} balance the conductivity
spectrum against two copies of the Dirichlet spectrum. Equivalently, we have
\[
N_N^{<}(s)-N_D^{<}(s)
\geq
1+\bigl|N_{\mathrm c}^{<}(s)-N_D^{<}(s)\bigr|.
\]
At the first Dirichlet threshold, the conductivity spectrum can be eliminated
by sharp isoperimetric estimates of Bucur--Henrot and
Bucur--Durus--Oudet
\cite{BucurHenrot2019,BucurDurusOudet2008}. A separate argument excludes the
remaining equality case and yields the strict inequality in
\cref{thm:planar-first}.

\subsection{Proof strategy}

The proof starts from two established ingredients: the Helmholtz--Hodge
decomposition on Lipschitz domains and the vector-field form involving
divergence and curl used in
\cite{Rohleder2025,FriesGoffengMiranda2026,Muravyev2026,HuaMunchZhang2026}.
The new point is to reorganize, rather than keep separate, the
finite-dimensional harmonic sector created by the holes. We absorb this
sector into the gradient space of functions whose traces are constant on each boundary component, with the
constants allowed to vary independently. This identifies the appropriate conductivity spectrum and leads to a relative Helmholtz decomposition valid for arbitrary planar
connectivity.

From the viewpoint of the existing Hodge-theoretic proofs, the obstruction on
multiply connected domains is therefore a spectral-bookkeeping issue rather
than a failure of the div--curl variational identity: the harmonic one-form
sector appears as a topological defect in the counting formula. The key
observation is that these degrees of freedom need not be counted as a separate
sector. Enlarging the scalar gradient sector to
$\nabla\mathcal V_{\mathrm c,\diamond}(\Omega)$ absorbs them into the
conductivity spectrum. In this way, topology is transferred from an index loss
to an auxiliary scalar spectrum, which can then be balanced against the
Dirichlet and Neumann spectra through the first-order factorization.

We then construct a first-order operator whose squared modulus on the scalar
side is the direct sum of the conductivity operator and the positive Neumann
Laplacian. In the opposite order, its quadratic form is the div--curl form on a
space of vector fields that contains two copies of every Dirichlet trial
space. The min--max principle gives \eqref{eq:counting-intro}. To deduce
\eqref{eq:planar-first}, we divide according to the relative positions of
$\gamma_1$ and $\lambda_1$. When $\gamma_1\leq\lambda_1$, sharp
isoperimetric inequalities and the strict ordering of the relevant Bessel
zeros give the result. When $\lambda_1<\gamma_1$, the counting theorem first
gives $\mu_3\leq\lambda_1$; equality is then ruled out by an
isotropic-subspace argument and the explicit Neumann component of the
first-order factorization.

Thus the classical ingredients are the conductivity problem, the Helmholtz--Hodge decomposition, the div--curl identity, the min--max principle,
and the sharp isoperimetric estimates. The new contributions are the
first-order factorization connecting the three spectra on multiply connected
domains, the counting inequality in \cref{thm:counting-intro}, and the
equality mechanism that gives \cref{thm:planar-first} without a simple-connectivity assumption.

\subsection{Paper organization}

In \cref{sec:conductivity} we define the conductivity spectrum and establish
the relative Helmholtz decomposition needed for multiply connected planar
domains. In \cref{sec:first-order} we construct the first-order operator and
prove \cref{thm:counting-intro} and \cref{cor:counting-gap}. The strict first-threshold inequality
\cref{thm:planar-first} is proved in \cref{sec:first-threshold}.

\section{Conductivity spectrum and a relative Helmholtz decomposition}
\label{sec:conductivity}

Throughout \cref{sec:conductivity,sec:first-order,sec:first-threshold},
$\Omega\subset\R^2$ is a bounded connected Lipschitz domain with boundary
components $\Gamma_0,\ldots,\Gamma_g$. We work in complex Hilbert spaces unless
we explicitly pass to real scalars. Set
\[
L^2_\diamond(\Omega):=
\left\{u\in L^2(\Omega):\int_\Omega u\dd x=0\right\},
\qquad
H^1_\diamond(\Omega):=H^1(\Omega)\cap L^2_\diamond(\Omega).
\]
Define
\[
\mathcal V_{\mathrm c}(\Omega):=
\left\{u\in H^1(\Omega):
\Tr u|_{\Gamma_\ell}\text{ is constant for }0\leq\ell\leq g\right\}
\]
and
$\mathcal V_{\mathrm c,\diamond}(\Omega):=
\mathcal V_{\mathrm c}(\Omega)\cap L^2_\diamond(\Omega)$.
The latter is a closed subspace of $H^1(\Omega)$, and its embedding into
$L^2_\diamond(\Omega)$ is compact. The densely defined closed positive form
in $L^2_\diamond(\Omega)$
\begin{equation}\label{eq:conductivity-form}
\mathfrak a_{\mathrm c}[u,v]
:=
\int_\Omega\nabla u\cdot\nabla\overline v\dd x,
\qquad
\dom\mathfrak a_{\mathrm c}:=
\mathcal V_{\mathrm c,\diamond}(\Omega),
\end{equation}
defines a positive self-adjoint operator $A_{\mathrm c}$ with compact
resolvent. We denote its eigenvalues by
$0<\gamma_1\leq\gamma_2\leq\cdots$.
This is the locally constant conductivity problem of
\cite{BucurDurusOudet2008}. On a smooth domain, an eigenfunction is constant on
each $\Gamma_\ell$ and has zero total normal flux through each boundary
component. Indeed, admissible test functions may have independently prescribed
boundary constants, while their mean can be corrected by an element of
$H_0^1(\Omega)$. The weak equation therefore gives the flux conditions;
integrating the interior equation eliminates the constant Lagrange multiplier
created by working in $L^2_\diamond(\Omega)$.
For a vector field $U=(U_1,U_2)$, write
$\curl U:=\partial_1U_2-\partial_2U_1$. We use the standard spaces
\[
H(\operatorname{div};\Omega)
:=\{U\in L^2(\Omega;\C^2):\operatorname{div}U\in L^2(\Omega)\}
\]
and
\[
H(\curl;\Omega)
:=\{U\in L^2(\Omega;\C^2):\curl U\in L^2(\Omega)\}.
\]
Their normal and tangential traces are understood weakly.

Let $\mathsf J:\C^2\to\C^2$ be the counterclockwise rotation
$\mathsf J(z_1,z_2):=(-z_2,z_1)$. If $\nu$ is the outer unit normal, set
$\tau:=\mathsf J\nu$. The next decomposition is the relative planar Helmholtz
decomposition in the form suited to \eqref{eq:conductivity-form}.

\begin{lemma}\label{lem:relative-helmholtz}
One has the orthogonal decomposition
\begin{equation}\label{eq:relative-helmholtz}
L^2(\Omega;\C^2)
=
\nabla\mathcal V_{\mathrm c,\diamond}(\Omega)
\mathbin{\widehat\oplus}
\mathsf J\nabla H^1_\diamond(\Omega),
\end{equation}
where $\widehat\oplus$ denotes the orthogonal direct sum in
$L^2(\Omega;\C^2)$.
\end{lemma}

\begin{proof}
Take $u\in\mathcal V_{\mathrm c,\diamond}(\Omega)$ and
$v\in H^1_\diamond(\Omega)$. Since $\mathsf J\nabla v\in H(\operatorname{div};\Omega)$
and $\operatorname{div}(\mathsf J\nabla v)=0$, the weak Green formula gives
\[
\int_\Omega\nabla u\cdot \mathsf J\nabla\overline v\dd x
=
\sum_{\ell=0}^g c_\ell
\left\langle \mathsf J\nabla\overline v\cdot\nu,1\right\rangle_{\Gamma_\ell},
\]
where $c_\ell$ is the trace of $u$ on $\Gamma_\ell$. The normal trace of
$\mathsf J\nabla v$ is, up to the choice of orientation, the tangential derivative of
$\Tr v$. Its pairing with the constant function on each closed boundary
component is zero. Hence the two spaces in \eqref{eq:relative-helmholtz} are
orthogonal.

We also record the closedness of the two summands. Suppose that
$\nabla u_n$ is a Cauchy sequence in $L^2(\Omega;\C^2)$ with
$u_n\in\mathcal V_{\mathrm c,\diamond}(\Omega)$. Since
$u_n-u_m$ has mean zero, the Poincar\'e inequality gives
\[
\|u_n-u_m\|_{H^1(\Omega)}
\leq C\|\nabla u_n-\nabla u_m\|_{L^2(\Omega)}.
\]
Thus $u_n\to u$ in $H^1(\Omega)$ for some $u\in H^1_\diamond(\Omega)$.
The trace map is continuous, and for every $\ell$ the one-dimensional
subspace of constant functions in $H^{1/2}(\Gamma_\ell)$ is closed.
Consequently, $\Tr u|_{\Gamma_\ell}$ is constant for every $\ell$, so
$u\in\mathcal V_{\mathrm c,\diamond}(\Omega)$. It follows that
$\nabla\mathcal V_{\mathrm c,\diamond}(\Omega)$ is closed in
$L^2(\Omega;\C^2)$. The same argument, using that $\mathsf J$ is unitary,
shows that $\mathsf J\nabla H^1_\diamond(\Omega)$ is closed.

It remains to identify the orthogonal complement of
$\mathsf J\nabla H^1_\diamond(\Omega)$. Let
$U\in L^2(\Omega;\C^2)$ be orthogonal to that space. Since subtracting a
constant does not change a gradient, $U$ is in fact orthogonal to
$\mathsf J\nabla H^1(\Omega)$. Testing with compactly supported functions gives
$\curl U=0$ in distributions. Hence $U\in H(\curl;\Omega)$, and the weak Green
formula, together with the surjectivity of the trace map, gives
$U\cdot\tau=0$ on $\partial\Omega$.

Choose a ball $B$ with $\overline\Omega\subset B$, and let $\widetilde U$ be
the extension of $U$ by zero to $B$. The zero-extension formula for
$H(\curl)$, or directly the weak Green formula, gives
\[
\curl\widetilde U=0
\qquad\text{in }\mathcal D'(B),
\]
because both $\curl U$ and the tangential trace of $U$ vanish. Since $B$ is a ball, the distributional Poincar\'e lemma yields an
$h\in H^1(B)$ such that $\widetilde U=\nabla h$. On every connected
component of $B\setminus\overline\Omega$ one has $\nabla h=0$, so $h$ is
constant there. As $h$ is a single $H^1(B)$ function, its traces from the two
sides of each boundary component agree. It follows that the trace of
$h|_\Omega$ is constant on each $\Gamma_\ell$. Thus
$h|_\Omega\in\mathcal V_{\mathrm c}(\Omega)$, and subtracting its mean shows
that $U\in\nabla\mathcal V_{\mathrm c,\diamond}(\Omega)$.
Together with the orthogonality already proved, this gives
\[
\bigl(\mathsf J\nabla H^1_\diamond(\Omega)\bigr)^\perp
=
\nabla\mathcal V_{\mathrm c,\diamond}(\Omega),
\]
which proves \eqref{eq:relative-helmholtz}.
\end{proof}

\section{First-order factorization and the three-spectrum counting inequality}
\label{sec:first-order}

Define the closed densely defined operator
\begin{equation}\label{eq:first-order-operator}
\begin{aligned}
\mathscr D:
L^2_\diamond(\Omega)\oplus L^2_\diamond(\Omega)
&\longrightarrow L^2(\Omega;\C^2),\\
\mathscr D(u,v)&:=\nabla u+\mathsf J\nabla v,
\end{aligned}
\qquad
\dom\mathscr D:=
\mathcal V_{\mathrm c,\diamond}(\Omega)
\oplus H^1_\diamond(\Omega).
\end{equation}
By \cref{lem:relative-helmholtz}, $\mathscr D$ is bijective and
\[
\|\mathscr D(u,v)\|_{L^2}^2
=
\|\nabla u\|_{L^2}^2+\|\nabla v\|_{L^2}^2.
\]

The following proposition makes the spectral content of $\mathscr D$ explicit.
Set
\[
\mathcal X_{\mathrm c}(\Omega):=
\left\{
\begin{array}{l|l}
U\in H(\operatorname{div};\Omega)\cap H(\curl;\Omega)
& U\cdot\tau=0\text{ on }\partial\Omega,\\[-1mm]
& \left\langle U\cdot\nu,1\right\rangle_{\Gamma_\ell}=0
\text{ for }0\leq\ell\leq g
\end{array}
\right\}.
\]

\begin{proposition}\label{prop:first-order-spectrum}
The adjoint of $\mathscr D$ is given by
\begin{equation}\label{eq:first-order-adjoint}
\dom\mathscr D^*=\mathcal X_{\mathrm c}(\Omega),
\qquad
\mathscr D^*U=
\bigl(-\operatorname{div}U,-\curl U\bigr).
\end{equation}
Moreover,
\[
\mathscr D^*\mathscr D=A_{\mathrm c}\oplus A_N^+,
\]
where $A_N^+$ is the Neumann Laplacian restricted to
$L^2_\diamond(\Omega)$. The operator
$\mathscr B:=\mathscr D\mathscr D^*$ is positive, has compact resolvent and
trivial kernel, and is associated with the closed form
\begin{equation}\label{eq:vector-field-form}
\mathfrak b[U,V]
:=
\int_\Omega
\left(
\operatorname{div}U\,\overline{\operatorname{div}V}
+
\curl U\,\overline{\curl V}
\right)\dd x,
\qquad
\dom\mathfrak b:=\mathcal X_{\mathrm c}(\Omega).
\end{equation}
Its eigenvalues, including multiplicities, are precisely
\begin{equation}\label{eq:union-spectrum}
\{\gamma_j:j\geq1\}
\mathbin{\uplus}
\{\mu_j:j\geq2\}.
\end{equation}
Here $\uplus$ denotes multiset union.
\end{proposition}

\begin{proof}
Let $U\in\dom\mathscr D^*$, and write
$\mathscr D^*U=(f,q)$ with $f,q\in L^2_\diamond(\Omega)$. For
$\varphi\in C_c^\infty(\Omega)$, set
\[
\varphi_\diamond
:=
\varphi-\frac{1}{|\Omega|}\int_\Omega\varphi\dd x.
\]
Then $\varphi_\diamond$ belongs both to
$\mathcal V_{\mathrm c,\diamond}(\Omega)$ and to
$H^1_\diamond(\Omega)$, and its gradient equals $\nabla\varphi$. Since $f$ and
$q$ have mean zero, the adjoint identity gives
\[
\int_\Omega\nabla\varphi\cdot\overline U\dd x
=
\int_\Omega\varphi\,\overline f\dd x,
\qquad
\int_\Omega\mathsf J\nabla\varphi\cdot\overline U\dd x
=
\int_\Omega\varphi\,\overline q\dd x.
\]
Consequently,
$\operatorname{div}U,\curl U\in L^2(\Omega)$ and
$f=-\operatorname{div}U$, $q=-\curl U$.

The weak Green formulas now read
\[
\int_\Omega\nabla u\cdot\overline U\dd x
=
-\int_\Omega u\,\overline{\operatorname{div}U}\dd x
+
\sum_{\ell=0}^g c_\ell
\left\langle\overline{U\cdot\nu},1\right\rangle_{\Gamma_\ell}
\]
for $u\in\mathcal V_{\mathrm c,\diamond}(\Omega)$, and
\[
\int_\Omega \mathsf J\nabla v\cdot\overline U\dd x
=
-\int_\Omega v\,\overline{\curl U}\dd x
+
\left\langle\Tr v,\overline{U\cdot\tau}\right\rangle_{\partial\Omega}
\]
for $v\in H^1_\diamond(\Omega)$. Fix
$\eta\in C_c^\infty(\Omega)$ with $\int_\Omega\eta\dd x\neq0$. After choosing
a function with prescribed boundary constants, or with an arbitrary prescribed
trace, one can subtract a suitable multiple of $\eta$ to make its mean zero
without changing its trace. Hence the boundary terms vanish for all admissible
$u$ and $v$ if and only if the conditions defining
$\mathcal X_{\mathrm c}(\Omega)$ hold.

Conversely, if $U\in\mathcal X_{\mathrm c}(\Omega)$, then the weak divergence
and Stokes formulas give
\[
\int_\Omega\operatorname{div}U\dd x
=
\sum_{\ell=0}^g\left\langle U\cdot\nu,1\right\rangle_{\Gamma_\ell}
=0,
\qquad
\int_\Omega\curl U\dd x
=
\left\langle U\cdot\tau,1\right\rangle_{\partial\Omega}
=0.
\]
Thus $(-\operatorname{div}U,-\curl U)$ belongs to
$L^2_\diamond(\Omega)\oplus L^2_\diamond(\Omega)$, and the two Green formulas
prove the adjoint identity. This establishes
\eqref{eq:first-order-adjoint}, including the converse implication.

The orthogonality identity following \eqref{eq:first-order-operator} shows, at
the level of closed forms, that
$\mathscr D^*\mathscr D=A_{\mathrm c}\oplus A_N^+$. In particular,
$\mathscr D^*\mathscr D$ has compact resolvent. Since $\mathscr D$ is injective
and onto, its polar decomposition has the form $
\mathscr D=\mathscr U(\mathscr D^*\mathscr D)^{1/2}$
with a unitary operator $\mathscr U$. Consequently, $
\mathscr B
=
\mathscr D\mathscr D^*
=
\mathscr U(\mathscr D^*\mathscr D)\mathscr U^*.
$
This standard closed-operator argument is also the elementary spectral
mechanism underlying Hilbert complexes; see
\cite{BruningLesch1992,FriesGoffengMiranda2026,Muravyev2026}. It proves
compactness of the resolvent, absence of a zero mode, and
\eqref{eq:union-spectrum}. Finally, the form of
$\mathscr D\mathscr D^*$ is exactly \eqref{eq:vector-field-form}.
\end{proof}

We now prove the counting inequality stated in the introduction.

\begin{proof}[Proof of \cref{thm:counting-intro}]
The form-domain inclusion
$\mathcal V_{\mathrm c,\diamond}(\Omega)
\subset H^1_\diamond(\Omega)$ and the min--max principle give
$\mu_{j+1}\leq\gamma_j$ for every $j\geq1$. Therefore
\begin{equation}\label{eq:first-counting-bound}
N_{N,+}^{<}(s)\geq N_{\mathrm c}^{<}(s).
\end{equation}

Let $E_D^{<}(s)$ be the span of the Dirichlet eigenfunctions with eigenvalue
strictly below $s$, and set
$\mathcal W_s:=E_D^{<}(s)\times E_D^{<}(s)$. Then
$\dim\mathcal W_s=2N_D^{<}(s)$ and
$\mathcal W_s\subset H_0^1(\Omega;\C^2)
\subset\mathcal X_{\mathrm c}(\Omega)$. For $U=(u_1,u_2)\in H_0^1(\Omega;\C^2)$, expanding the
squares and integrating the mixed terms by parts, first for
$U\in C_c^\infty(\Omega;\C^2)$ and then by density, gives the
planar div--curl identity
\begin{equation}\label{eq:gaffney}
\int_\Omega
\bigl(
|\operatorname{div}U|^2+|\curl U|^2
\bigr)\dd x
=
\int_\Omega
\bigl(
|\nabla u_1|^2+|\nabla u_2|^2
\bigr)\dd x;
\end{equation}
see also \cite[Equation~(16)]{Rohleder2025}.
The combination of the two-copy Dirichlet trial space and the div--curl
identity forms a standard variational device in Rohleder's argument and its
Hodge-theoretic reformulations and generalizations
\cite{Rohleder2025,FriesGoffengMiranda2026,Muravyev2026}. What is specific
here is that \eqref{eq:union-spectrum} replaces the usual
Dirichlet--Neumann--harmonic spectral decomposition by the
conductivity--Neumann spectrum.
Thus every nonzero $U\in\mathcal W_s$ satisfies
$\mathfrak b[U]<s\|U\|_{L^2}^2$. The min--max principle for $\mathscr B$ and
\eqref{eq:union-spectrum} yield
\begin{equation}\label{eq:second-counting-bound}
N_{\mathrm c}^{<}(s)+N_{N,+}^{<}(s)
\geq 2N_D^{<}(s).
\end{equation}
Combining \eqref{eq:first-counting-bound} and
\eqref{eq:second-counting-bound} proves \eqref{eq:counting-intro}.
\end{proof}

The same estimate has a useful symmetric formulation.

\begin{corollary}\label{cor:counting-gap}
For every $s>0$,
\begin{equation}\label{eq:counting-gap}
N_N^{<}(s)-N_D^{<}(s)
\geq
1+\bigl|N_{\mathrm c}^{<}(s)-N_D^{<}(s)\bigr|.
\end{equation}
\end{corollary}

\begin{proof}
Subtracting $N_D^{<}(s)$ from the two alternatives in
\eqref{eq:counting-intro} gives
\[
N_{N,+}^{<}(s)-N_D^{<}(s)
\geq
\bigl|N_{\mathrm c}^{<}(s)-N_D^{<}(s)\bigr|.
\]
Adding the Neumann zero mode proves \eqref{eq:counting-gap}.
\end{proof}

\section{The first Dirichlet threshold in the plane}
\label{sec:first-threshold}

We first isolate a standard elementary consequence of the spectral theorem
that will be used to rule out equality.

\begin{lemma}\label{lem:isotropic-subspace}
Let $T$ be a self-adjoint operator with compact resolvent and an associated
closed form $\mathfrak t$ that is bounded from below. Fix $a\in\R$, and assume that the spectral subspace
of $T$ in $(-\infty,a)$ is one-dimensional. If a two-dimensional subspace
$W\subset\dom\mathfrak t$ satisfies
\[
\mathfrak t[u,v]=a(u,v)
\qquad\text{for all }u,v\in W,
\]
then $W\cap\ker(T-a)\neq\{0\}$.
\end{lemma}

\begin{proof}
Let $P_-$ be the spectral projection of $T$ onto $(-\infty,a)$. Since
$\dim W=2$ and $\operatorname{rank}P_-=1$, there exists
$0\neq w\in W$ such that $P_-w=0$. The hypothesis gives
$\mathfrak t[w]=a\|w\|^2$. In the spectral representation of $T$, the form of
$T-a$ is nonnegative on $(I-P_-)\dom\mathfrak t$, and it vanishes only on
$\ker(T-a)$. Hence $w\in\ker(T-a)$.
\end{proof}

We can now prove the strict first-threshold inequality.

\begin{proof}[Proof of \cref{thm:planar-first}]
Let $J_0$ and $J_1$ denote the Bessel functions of the first kind of orders
$0$ and $1$, respectively. Let $j_{0,1}$ be the first positive zero of $J_0$,
and let $j'_{1,1}$ be the first positive zero of $J_1'$. We distinguish the
position of the first conductivity eigenvalue.

Assume first that $\gamma_1\leq\lambda_1$. We first verify that the
conductivity inequality of \cite{BucurDurusOudet2008} applies to the form
\eqref{eq:conductivity-form}. The conductivity space used there is the
$H^1(\Omega)$-closure of the restrictions of functions
$w\in H^1_{\mathrm{loc}}(\R^2)$ such that $\nabla w=0$ in a neighbourhood of
$\R^2\setminus\Omega$. For a bounded connected planar Lipschitz domain, this
closed space is exactly $\mathcal V_{\mathrm c}(\Omega)$. Indeed, every
function in the defining dense class has constant trace on each
$\Gamma_\ell$, and this property is closed under $H^1$ convergence by the
continuity of the trace map. Conversely, if
$u\in\mathcal V_{\mathrm c}(\Omega)$ has trace $c_\ell$ on
$\Gamma_\ell$, choose $h\in H^1_{\mathrm{loc}}(\R^2)$ which is equal to
$c_\ell$ on the corresponding complementary component and in a collar of
$\Gamma_\ell$. Then
\[
u-h|_\Omega\in H_0^1(\Omega),
\]
and approximating this difference by functions in
$C_c^\infty(\Omega)$ proves the reverse inclusion. Thus the zero-mean
conductivity form domain in \cite{BucurDurusOudet2008} is precisely
$\mathcal V_{\mathrm c,\diamond}(\Omega)$.

The sharp theorems of Bucur--Henrot and Bucur--Durus--Oudet therefore
give
\[
\mu_3(\Omega)|\Omega|
\leq 2\pi(j'_{1,1})^2,
\qquad
\gamma_1(\Omega)\geq\gamma_1(B_1\cup B_2),
\]
where $B_1$ and $B_2$ are two disjoint equal disks of total area
$|\Omega|$. Here $\gamma_1(B_1\cup B_2)$ denotes the first eigenvalue
of the Dirichlet energy form on the zero-mean conductivity space of
\cite{BucurDurusOudet2008} for this disconnected open set; see
\cite{BucurHenrot2019,BucurDurusOudet2008}. Bucur and Henrot enumerate the
Neumann spectrum as
\[
0=\mu_0\leq\mu_1\leq\mu_2\leq\cdots,
\]
so their $\mu_2$ is our $\mu_3$.

It remains to record the conductivity eigenvalue of $B_1\cup B_2$. Since its
complement is connected, the conductivity space reduces to the global
constant-boundary space studied in \cite{GrecoLucia2008}; an admissible
function has one common boundary value $c$ on the two circles. For clarity,
the short computation is included. After identifying both disks with a disk
$B$, write the two restrictions as $u_1,u_2$ and set
\[
u_+:=\frac{u_1+u_2}{2},
\qquad
u_-:=\frac{u_1-u_2}{2}.
\]
Then
\[
u_-\in H_0^1(B),
\qquad
u_+-c\in H_0^1(B),
\]
and the zero-mean condition gives
\[
\int_B u_+\dd x=0.
\]
Consequently,
\[
\|u_+-c\|_{L^2(B)}^2
=
\|u_+\|_{L^2(B)}^2+|B||c|^2
\geq
\|u_+\|_{L^2(B)}^2.
\]
The Dirichlet Poincar\'e inequality, together with the orthogonal splitting
of both the $L^2$ norm and the Dirichlet energy into the symmetric and
antisymmetric parts, yields
\[
\int_{B_1\cup B_2}|\nabla u|^2\dd x
\geq
\lambda_1(B)\int_{B_1\cup B_2}|u|^2\dd x.
\]
Equality is attained by the antisymmetric test
$(\phi_1,-\phi_1)$, where $\phi_1$ is a first Dirichlet eigenfunction of
$B$. Since $|B|=|\Omega|/2$, it follows that
\[
\gamma_1(B_1\cup B_2)
=
\lambda_1(B)
=
\frac{2\pi j_{0,1}^2}{|\Omega|}.
\]
Consequently,
\begin{equation}\label{eq:sharp-isoperimetric}
\mu_3(\Omega)|\Omega|
\leq 2\pi(j'_{1,1})^2,
\qquad
\gamma_1(\Omega)|\Omega|
\geq 2\pi j_{0,1}^2.
\end{equation}
The standard recurrence and zero-interlacing relations for Bessel functions
give
$j'_{1,1}<j_{0,1}$. Indeed, $j_{0,1}<j_{1,1}$, while
$J_1'(x)=J_0(x)-x^{-1}J_1(x)$. Hence
$J_1'(j_{0,1})<0$, whereas $J_1'(x)\to 1/2$ as $x\downarrow0$.
Thus \(j'_{1,1}<j_{0,1}\); see \cite[\S\S10.6(i), 10.7(i), and 10.21(i)]{DLMF}. Thus
\eqref{eq:sharp-isoperimetric} yields
\[
\mu_3(\Omega)<\gamma_1(\Omega)\leq\lambda_1(\Omega).
\]

It remains to consider $\lambda_1<\gamma_1$. Since $\lambda_1$ is simple,
choose
$\lambda_1<s<\min\{\lambda_2,\gamma_1\}$. Then
$N_D^{<}(s)=1$ and $N_{\mathrm c}^{<}(s)=0$, so
\eqref{eq:counting-intro} gives $N_{N,+}^{<}(s)\geq2$. Letting
$s\downarrow\lambda_1$ shows that $\mu_3\leq\lambda_1$.

Suppose, for a contradiction, that equality holds, and put $t:=\lambda_1=\mu_3$. Filonov's strict one-step inequality gives $\mu_2<t$ \cite{Filonov2005}. Since $t<\gamma_1$, the spectral description~\eqref{eq:union-spectrum} shows that $\mathscr B$ has exactly one eigenvalue below $t$. Let $\phi_1$ be a real first Dirichlet eigenfunction and, in the real Hilbert space, set \(W:=\{(a\phi_1,b\phi_1):a,b\in\R\}\). The polarized form of \eqref{eq:gaffney} gives $\mathfrak b[U,V]=t(U,V)_{L^2}$ for all $U,V\in W$. By~\cref{lem:isotropic-subspace}, there is a nonzero $U=(a\phi_1,b\phi_1)\in W\cap\ker(\mathscr B-t)$.

The eigenspace at $t$ comes only from the Neumann summand in
\eqref{eq:union-spectrum}. The unitary equivalence in the proof of
\cref{prop:first-order-spectrum} therefore gives
$U=\mathsf J\nabla\psi$ for a Neumann eigenfunction $\psi$ with eigenvalue $t$.
In particular, $\operatorname{div}U=0$, and hence
\[
a\,\partial_1\phi_1+b\,\partial_2\phi_1=0
\quad\text{in }\Omega,
\qquad (a,b)\neq(0,0).
\]
Let $\widetilde\phi_1$ be the zero extension of $\phi_1$ to $\R^2$.
Because $\phi_1\in H_0^1(\Omega)$, one has
$\widetilde\phi_1\in H^1(\R^2)$, and the same directional derivative vanishes
in distributions on $\R^2$. After a rotation,
$\partial_1\widetilde\phi_1=0$. Hence
$\widetilde\phi_1$ is independent of the first coordinate almost everywhere;
an $L^2(\R^2)$ function with this property must vanish. This contradicts
$\phi_1\neq0$, so equality is impossible.
\end{proof}

\section*{Statement on AI usage}

The overall strategy and the main mathematical ideas of this paper were
developed by the authors. In particular, the conductivity-spectrum
viewpoint, the relative Helmholtz decomposition, and the first-order
framework leading to the three-spectrum comparison were conceived by the
authors.

ChatGPT was used to assist with parts of the proofs. Most notably,
a substantial part of the proof of
Proposition~\ref{prop:first-order-spectrum} was initially developed with
the assistance of ChatGPT, starting from the first-order operator and
spectral framework proposed by the authors. ChatGPT was also used to help
formalize and check some technical steps in
Lemma~\ref{lem:relative-helmholtz} and in the equality analysis in the
proof of Theorem~\ref{thm:planar-first}. The resulting arguments were
subsequently revised and checked carefully by the authors. The authors
take full responsibility for all content of the paper.


\begin{thebibliography}{99}
\small
\raggedright

\bibitem{ArendtMazzeo2012}
W.~Arendt and R.~Mazzeo,
\emph{Friedlander's eigenvalue inequalities and the Dirichlet-to-Neumann semigroup},
Commun. Pure Appl. Anal. \textbf{11} (2012), no.~6, 2201--2212.

\bibitem{Aviles1986}
P.~Aviles,
\emph{Symmetry theorems related to Pompeiu's problem},
Amer. J. Math. \textbf{108} (1986), no.~5, 1023--1036.

\bibitem{BenguriaLevitinParnovski2009}
R.~Benguria, M.~Levitin, and L.~Parnovski,
\emph{Fourier transform, null variety, and Laplacian's eigenvalues},
J. Funct. Anal. \textbf{257} (2009), no.~7, 2088--2123.

\bibitem{BruningLesch1992}
J.~Br\"uning and M.~Lesch,
\emph{Hilbert complexes},
J. Funct. Anal. \textbf{108} (1992), no.~1, 88--132.

\bibitem{BucurDurusOudet2008}
D.~Bucur, I.~Durus, and E.~Oudet,
\emph{The eigenvalue problem with conductivity boundary conditions},
Control Cybernet. \textbf{37} (2008), no.~4, 811--829.

\bibitem{BucurHenrot2019}
D.~Bucur and A.~Henrot,
\emph{Maximization of the second non-trivial Neumann eigenvalue},
Acta Math. \textbf{222} (2019), no.~2, 337--361.

\bibitem{CoxMacLachlanSteeves2019}
G.~Cox, S.~MacLachlan, and L.~Steeves,
\emph{Isoperimetric relations between Dirichlet and Neumann eigenvalues},
preprint, arXiv:1906.10061 (2019).

\bibitem{Filonov2005}
N.~Filonov,
\emph{On an inequality between Dirichlet and Neumann eigenvalues for the
Laplace operator},
Algebra i Analiz \textbf{16} (2004), no.~2, 172--176;
English transl., St. Petersburg Math. J. \textbf{16} (2005), no.~2, 413--416.

\bibitem{Filonov2026}
N.~Filonov,
\emph{Inequalities between Dirichlet and Neumann eigenvalues in large
dimensions},
preprint, arXiv:2606.30120 (2026).

\bibitem{FreitasGama}
P.~Freitas and M.~Gama,
\emph{On the (growing) gap between Dirichlet and Neumann eigenvalues},
Indiana Univ. Math. J., to appear; arXiv:2405.18079v2 (2024).

\bibitem{Friedlander1991}
L.~Friedlander,
\emph{Some inequalities between Dirichlet and Neumann eigenvalues},
Arch. Ration. Mech. Anal. \textbf{116} (1991), no.~2, 153--160.

\bibitem{FriesGoffengMiranda2026}
M.~Fries, M.~Goffeng, and G.~Miranda,
\emph{Reproving Friedlander's inequality with the de Rham complex},
J. Spectr. Theory \textbf{15} (2025), no.~4, 1593--1613;
corrected version under the title
\emph{Friedlander's inequality and the de Rham complex},
arXiv:2412.03369v4, 2026.

\bibitem{GesztesyMitrea2009}
F.~Gesztesy and M.~Mitrea,
\emph{Nonlocal Robin Laplacians and some remarks on a paper by Filonov on
eigenvalue inequalities},
J. Differential Equations \textbf{247} (2009), no.~10, 2871--2896.

\bibitem{GrecoLucia2008}
A.~Greco and M.~Lucia,
\emph{Laplacian eigenvalues for mean zero functions with constant Dirichlet
data},
Forum Math. \textbf{20} (2008), no.~5, 763--782.

\bibitem{Hatcher2025}
L.~Hatcher,
\emph{Geometric inequalities between Dirichlet and Neumann eigenvalues},
preprint, arXiv:2504.18517 (2025).

\bibitem{HuaMunchZhang2026}
B.~Hua, F.~M\"unch, and H.~Zhang,
\emph{Inequalities between Dirichlet and Neumann eigenvalues on surfaces},
Math. Ann. \textbf{395} (2026), article~49.

\bibitem{LevineWeinberger1986}
H.~A. Levine and H.~F. Weinberger,
\emph{Inequalities between Dirichlet and Neumann eigenvalues},
Arch. Ration. Mech. Anal. \textbf{94} (1986), no.~3, 193--208.

\bibitem{LevitinMangoubiPolterovich2023}
M.~Levitin, D.~Mangoubi, and I.~Polterovich,
\emph{Topics in Spectral Geometry},
Graduate Studies in Mathematics, vol.~237,
American Mathematical Society, Providence, RI, 2023.

\bibitem{Muravyev2026}
M.~Muravyev,
\emph{Hodge Laplacian eigenvalues on surfaces with boundary},
Ann. Math. Qu\'ebec \textbf{50} (2026), 73--92.

\bibitem{DLMF}
NIST Digital Library of Mathematical Functions,
Release 1.2.7 of 2026-06-15,
F.~W.~J. Olver, A.~B. Olde Daalhuis, D.~W. Lozier,
B.~I. Schneider, R.~F. Boisvert, C.~W. Clark, B.~R. Miller,
B.~V. Saunders, H.~S. Cohl, and M.~A. McClain, eds.,
\url{https://dlmf.nist.gov/}.

\bibitem{Payne1955}
L.~E. Payne,
\emph{Inequalities for eigenvalues of membranes and plates},
J. Ration. Mech. Anal. \textbf{4} (1955), 517--529.

\bibitem{Polya1952}
G.~P\'olya,
\emph{Remarks on the foregoing paper},
J. Math. Phys. \textbf{31} (1952), 55--57.

\bibitem{Rohleder2025}
J.~Rohleder,
\emph{Inequalities between Neumann and Dirichlet Laplacian eigenvalues on
planar domains},
Math. Ann. \textbf{392} (2025), no.~4, 5553--5571.

\bibitem{Safarov2008}
Y.~Safarov,
\emph{On the comparison of the Dirichlet and Neumann counting functions},
in \emph{Spectral Theory of Differential Operators: M.~Sh.~Birman 80th
Anniversary Collection}, Amer. Math. Soc. Transl. Ser.~2, vol.~225,
American Mathematical Society, Providence, RI, 2008, 191--204.

\bibitem{SafarovFilonov2010}
Y.~G. Safarov and N.~D. Filonov,
\emph{Asymptotic estimates of the difference between the Dirichlet and Neumann counting functions},
Funct. Anal. Appl. \textbf{44} (2010), no.~4, 286--294.

\bibitem{Schmatzler2026}
T.~Schmatzler,
\emph{The Levine--Weinberger and Friedlander--Filonov inequalities for some
classes of elliptic operators},
preprint, arXiv:2506.09468v2 (2026).

\end{thebibliography}
\end{document}